\documentclass{article}

\usepackage{subfiles} 

\usepackage{packages/commonpackage}
\usepackage{packages/package_for_list_of_content}
\usepackage{packages/package_to_delete}
\usepackage{packages/macro}

\numberwithin{equation}{section}
\theoremstyle{plain}
\newtheorem{theorem}{Theorem}[section]
\newtheorem{proposition}[theorem]{Proposition}

\newtheorem{lemma}[theorem]{Lemma}
\newtheorem{definition}[theorem]{Definition}
\newtheorem{remark}[theorem]{Remark}
\newtheorem{example}[theorem]{Example}
\newtheorem{claim}[theorem]{Claim}

\newtheorem{conjecture}[theorem]{Conjecture}

\title{Rational points near planar flat curves}
\author{Mingfeng Chen}
\date{}

\begin{document}
\maketitle

\begin{abstract}
    We establish asymptotic formulas for counting rational points near finite type curves on $\R^2$, generalizing Huang's result \cite{Hua15}.
\end{abstract}

\section{Introduction}
\subsection{Statement of the main theorem}
Let $f:I \to \R$ be a real-valued function, where $I$ is a compact interval in $\R$. Given $Q\ge 1$ as an integer and $\delta \in (0,\frac{1}{2})$, consider
 \begin{align}
    N_f(Q,\delta):=\# \{(a,q)\in \Z \times \N : 1\leq q\leq Q,\ \frac{a}{q}\in I,\ ||qf(\frac{a}{q})||<\delta   \},
\end{align}
where $||x||:=\min_{k\in \Z} |x-k|$. The quantity $N_f(Q,\delta)$ counts the number of rational points near the planar curve $(t,f(t))$. 
\begin{definition}
    Let $d\ge 2$ be an integer and $C(t)=(t,f(t))$ be a curve defined near a small neighbourhood of $a_0$. We say that $C$ is of type $d$ at $a_0$ if 
    \begin{align}
        f^{(k)}(a_0)=0
    \end{align}
    for $1\leq k<d$ and 
    \begin{align}
        f^{(d)}(a_0)\neq 0.
    \end{align}
\end{definition}
Typical examples of finite type curves include $f(x)=x^d$ and the Fermat curve discussed below. Our main result is the following asymptotic formula for counting rational points near finite type curves:
\begin{theorem}\label{Thm1.2}
    Let $C(t)=(t,f(t))$ be a $C^{d+1}$ curve of type $d\ge 3$ at $a_0$, then there exists $\epsilon_0>0$, such that for any interval $I\subseteq (a_0-\epsilon_0,a_0+\epsilon_0)$, we have
    \begin{align}\label{1.12}
         N_f(Q,\delta)=|I|\delta Q^2+O\Big{(}\delta^{\frac{1}{d}}Q^{2-\frac{1}{d}}+Q^{1+\epsilon }\Big{)}
    \end{align}
    for all integers $Q\ge 1$ and $\delta \in (Q^{-1},\frac{1}{2})$. In other words, we have
    \begin{align}
        N_f(Q,\delta)=(|I|+o(1))\delta Q^2
    \end{align}
    when $\delta \in (Q^{-\frac{1}{d-1}+\epsilon},\frac{1}{2})$ for some $\epsilon>0$ and $Q \to \infty$.
\end{theorem}
\begin{remark}
    When $d=2$, our proof recovers the main result of Huang in \cite{Hua15}, which is 
    \begin{align}
         N_f(Q,\delta)=|I|\delta Q^2+O\Big{(}\delta^{\frac{1}{2}}(\log\frac{1}{\delta}) Q^{\frac{3}{2}}+Q^{1+\epsilon }\Big{)}
    \end{align}
    for strictly convex $f\in C^3(I)$. Our result is a generalization of Huang's result to finite type curves on $\R^2$.
\end{remark}
Throughout this paper, we use the notation $A\lesim B$ to indicate that there exists a constant $C>0$, such that $A\leq CB$. Similarly, $A\simeq B$ means $A\lesim B$ and $B\lesim A$. The implicit constant depends only on the underlying manifold and some fixed parameters, and it is in particular independent of $Q$ and $\delta$.

Let us clarify the terms that appear on the right hand side of \eqref{1.12}.
\begin{example}\label{eg1.3}
    Consider the parabola $g(x)=x^2,x\in [0,1]$. The rational points $(\frac{a}{q},\frac{a^2}{q^2})$, with $1\leq a \leq q$ and $1\leq q\leq \sqrt{Q}$ are on the curve $(x,g(x))$. So we have
    \begin{align}
        N_{g}(Q,\delta)\ge N_g(Q,0) \ge \sum_{1\leq q\leq \sqrt{Q}} q \gtrsim Q .
    \end{align}\end{example}
This example suggests that when $d=2$, the constant term $Q^{1+\epsilon}$ is optimal. For general $d\ge 3$, determining the best constant term in $Q$ turns out to be a different problem, which is not the main focus of this paper.
\begin{example}[\cite{Hua20}, Example 4]\label{eg1.4}
    Consider the Fermat curve $\mathcal{F}_d$:
    \begin{align}
        x^d+y^d=1.
    \end{align}
Then $\mathcal{F}_d$ has contact order $d-1$ with the lines $x=1$ and $y=1$, which gives
\begin{align}
    N_{\mathcal{F}_d}(Q,\delta) \gtrsim \delta^{\frac{1}{d}}Q^{2-\frac{1}{d}}.
\end{align}
\end{example}
This example suggests that the middle term $\delta^{\frac{1}{d}}Q^{2-\frac{1}{d}}$ is optimal.

\subsection{Preliminary results}
In \cite{Hux94}, Huxley employed the determinant method and some geometric arguments to obtain a nearly optimal upper bound for $C^2$ curves with a second-order derivative bounded below:
\begin{align}\label{Hux}
    N_f(Q,\delta) \lesim_{\epsilon} \delta^{1-\epsilon}Q^2+Q\log Q .
\end{align}
Later, Vaughan and Velani \cite{VV06} applied methods from analytic number theory to establish a sharp upper bound for planar $C^3$ strictly convex curves
\begin{align}
    N_{f}(Q,\delta)\lesim_{\epsilon} \delta Q^2+Q^{1+\epsilon} .
\end{align}
Notably, Huxley's bound \eqref{Hux} played a crucial role in the proof presented in \cite{VV06}.

In \cite{Hua15}, based on \cite{Hux94} and \cite{VV06}, Huang proved the following asymptotic formula:
\begin{theorem}[\cite{Hua15}, Theorem 3]
    Let $f:I \to \R$ be a $C^2$ function with Lipschitz continuous second derivative, assume
    \begin{align}
        c_1\leq |f''(x)|\leq c_2
    \end{align}
    for all $x\in I$. Then, for any $\epsilon>0$, $\delta\in (0,\frac{1}{2})$ and integer $Q>1$, we have 
    \begin{align}
        N_f(Q,\delta)=|I|\delta Q^2+O\Big{(} \delta^{\frac{1}{2}}(\log \frac{1}{\delta})Q^{\frac{3}{2}}+Q^{1+\epsilon} \Big{)} ,
    \end{align}
    where the implicit constant depends only on $I,c_1,c_2,\epsilon$ and the Lipschitz constant, and is in particular independent of $f,\delta$ and $Q$. In other words, we have
    \begin{align}
        N_f(Q,\delta)=(|I|+o(1))\delta Q^2
    \end{align}
    when $\delta \in (Q^{-1+\epsilon},\frac{1}{2})$ for some $\epsilon>0$ and $Q \to \infty$.
\end{theorem}

\subsection{More results in higher dimensions}
Let $\mathcal{M}$ be a bounded submanifold of $\R^n$ of dimension $m$. Given $Q\ge 1$ as an integer and $\delta\in (0,\frac{1}{2})$, consider
\begin{align}
    N_{\mathcal{M}}(Q,\delta):=\# \{ \frac{\mathbf{p}}{q}\in \Q^n: 1\leq q\leq Q, \dist_{\infty}(\frac{\mathbf{p}}{q},\mathcal{M})\leq \frac{\delta}{q} \} ,
\end{align}
where $\dist_{\infty}(\cdot,\cdot)$ represents the $L^{\infty}$-distance on $\R^n$. 

A trivial upper bound for $N_{\mathcal{M}}(Q,\delta)$ is 
\begin{align}
    N_{\mathcal{M}}(Q,\delta) \lesim Q^{m+1} .
\end{align}

\begin{example}
    Let $\mathcal{M}$ be a manifold that contains a rational $m$-linear subspace of $\R^n$, then we have
    \begin{align}
        N_{\mathcal{M}}(Q,\delta)\simeq Q^{m+1} .
    \end{align}
\end{example}
This example demonstrates that the trivial bound can be achieved if no additional conditions are imposed on the underlying manifold. However, a probabilistic heuristic suggests that we should expect
\begin{align}
    N_{\mathcal{M}}(Q,\delta) \simeq \delta^k Q^{m+1} ,
\end{align}
where $k=n-m$ is the codimension of $\mathcal{M}$. In the breakthrough work \cite{Ber12}, Beresnevich established the following sharp lower bound 
\begin{align}
    N_{\mathcal{M}}(Q,\delta) \gtrsim \delta^k Q^{m+1}
\end{align}
for $\delta \in (Q^{-\frac{1}{k}},\frac{1}{2})$ and $\mathcal{M}$ being any analytic submanifold containing at least one nondegenerate point(see \cite{Ber12} for a precise definition).

For the upper bound, it is reasonable to believe that if the manifold is properly curved, the probabilistic heuristic bound should hold. The main conjecture is stated as follows(see \cite{Hua20} and references within)
\begin{conjecture}[\cite{Hua20}, Conjecture 1]
    For any compact submanifold of $\R^n$ with proper curvature conditions, we have
\begin{align}
    N_{\mathcal{M}}(Q,\delta) \simeq c_{\mathcal{M}}\delta^k Q^{m+1}
\end{align}
when $\delta \in (Q^{-\frac{1}{k}+\epsilon},\frac{1}{2})$ for some $\epsilon>0$ and $Q \to \infty$, where $m$ is the dimension of the manifold and $k=n-m$ is the codimension.
\end{conjecture}
Such estimates have numerous applications to other problems in number theory, including Diophantine inequalities, Serre's dimension growth conjecture and metric Diophantine approximation on manifolds(see \cite{Hua20} and \cite{ST23}). 

For manifolds with dimension $m\ge 2$, numerous results have been established. In Huang's groundbreaking work \cite{Hua20}, a sharp bound was proven for rational points near hypersurfaces with non-vanishing Gaussian curvature:
\begin{theorem} [\cite{Hua20}, Theorem 1]
    Let $n\ge 3$, $\mathcal{M}$ be a compact $C^l$ hypersurface in $\R^n$ with Gaussian curvature bounded away from zero. Then
    \begin{align}
        N_{\mathcal{M}}(Q,\delta)   \lesim \delta Q^n+E_n(Q),
    \end{align}
    where 
    \begin{align} 
        &l=\max \Big{(} \Big{\lfloor} \frac{n-1}{2}\Big{\rfloor}+5 ,n+1 \Big{)} ,  \\
        &E_3(Q):=Q^2 \exp(c\sqrt{\log Q}) ,        \\
        &E_n(Q):=Q^{n-1}(\log Q)^{\kappa}, \ n\ge 4.
    \end{align}
    $c,\kappa$ and the implicit constant only depend on $\mathcal{M}$.
\end{theorem}
Later, Schindler and Yamagishi \cite{SY22} generalized Huang's result to manifolds of high codimension under certain natural curvature conditions. The results of Huang \cite{Hua20}, as well as those of Schindler and Yamagishi \cite{SY22}, can be viewed as  non-degenerate cases, where the underlying manifold is not flat at every point. This allows the application of the stationary phase principle to the associated oscillatory integral.

Recently, Srivastava and Technau \cite{ST23} established sharp estimates for certain degenerate homogeneous hypersurfaces, provided their Gaussian curvature is non-vanishing away from the origin.
 
\begin{definition}
    For a real number $d\neq 0$, we denote $\mathcal{H}_d^0(\R^{n-1})$ to be the set of all smooth functions $f:\R^{n-1} \to \R$ which are homogeneous of degree $d$ and Hessian remains invertible away from the origin.
\end{definition}
The simplest example is $f(\bfx):=||\bfx||_2^d$. Their main result is the following:
\begin{theorem}[\cite{ST23}, Theorem 1.5]
    Let $n\ge 3$ be an integer and $d>\frac{2(n-1)}{2n-3}$ be a real number. Fix $\epsilon>0$ and $f\in \mathcal{H}_d^0(\R^{n-1})$. Then
    \begin{align}
        N_f(Q,\delta) \simeq \delta Q^n+\Big{(}\frac{\delta}{Q} \Big{)}^{\frac{n-1}{d}} Q^n
    \end{align}
    for any $Q\ge 1$ and $\delta \in (Q^{-1+\epsilon},\frac{1}{2})$, provided $\frac{2(n-1)}{2n-3}<d<n-1$.\\
    If $d\ge n-1$, then 
    \begin{align}
        \delta Q^n+\Big{(} \frac{\delta}{Q} \Big{)}^{\frac{n-1}{d}}Q^n \lesim N_f(Q,\delta) \lesim \delta Q^n+\Big{(} \frac{\delta}{Q} \Big{)}^{\frac{n-1}{d}}Q^{n+k\epsilon}    
        \end{align}
    for any $Q\ge 1$ and $\delta \in (Q^{-1+\epsilon},\frac{1}{2})$, with $k=k(n)>0$ depending only on n. The other implicit constants depend on $\epsilon,f,n$ and $d$, but not on $\delta$ and Q.
\end{theorem}

\subsection{Acknowledgments}
The author would like to thank Shaoming Guo for writing suggestions that improve the presentation of this paper. The author would like to thank Rajula Srivastava and Niclas Technau for helpful discussions.

\section{Applications}
In this section, we discuss applications of our main result Theorem 1.2, to Diophantine Approximation and a smooth variant of Serre's Dimension Growth Conjecture. We refer to \cite{Hua20} and \cite{ST23} for similar discussions and more details.
\subsection{Diophantine Approximation on Manifolds}
Given an approximation function $\psi: \N \to (0,\frac{1}{2})$ be decreasing, define the set of $\psi$-approximable numbers as
\begin{align}
    \mathcal{W}_n(\psi):=\{\bfx\in [0,1]^n:\max_{1\leq i\leq n}||qx_i||<\psi(q) \ \text{for infinitely many integers $q\ge 1$} \} .
\end{align}
Denoting by $m_n$ the n-dimensional Lebesgue measure on $[0,1]^n$, then we have the following 0-1 law:
\begin{theorem}[Khintchine's theorem; 1926]
    Let $\psi$ be as above, then 
    \begin{align}
        m_n(\mathcal{W}_n(\psi))=\left\{\begin{array}{lr}
        1  & \text{if}\ \sum_{q\ge 1}\psi(q)^n \ \text{diverges,} \\
        0 & \text{if}\ \sum_{q\ge 1}\psi(q)^n \ \text{converges.}    \end{array}
    \right.               \end{align}
\end{theorem}
We also have the following Hausdorff measure refinements of Khintchine's theorem:
\begin{theorem}[Jarn\'{i}k’s theorem; 1931]
    For any monotonic function $\psi$ and $0<s<n$,
\begin{align}
\mathcal{H}^s(\mathcal{W}_n(\psi))=\left\{\begin{array}{lr}
        \infty  & \text{if}\ \sum_{q\ge 1}q^n(\frac{\psi(q)}{q})^s \ \text{diverges,} \\
        0 & \text{if}\ \sum_{q\ge 1}q^n(\frac{\psi(q)}{q})^s\ \text{converges.}    \end{array}
    \right.     \end{align} where $\mathcal{H}^s$ denotes the s-dimensional Hausdorff measure. 
\end{theorem}
Our Theorem 1.2 gives the following result:
\begin{theorem}
    Let $s>\frac{1}{2}$, $\mathcal{C}$ be a $C^{d+1}$ curve of type $d\ge 2$, and $\epsilon>0$. Then, for any approximation function $\psi$, we have
    \begin{align}
        \mathcal{H}^s(\mathcal{W}_n(\psi) \cap \mathcal{C})=0 \ \text{if}\ \sum_{q\ge 1} (\frac{\psi(q)}{q})^{s+1}q^2 <\infty \ \text{and} \ \sum_{q\ge 1}(\frac{\psi(q)}{q})^{s+\frac{1}{d}+\epsilon}q<\infty.
    \end{align}
\end{theorem}

\subsection{A smooth variant of Serre's Dimension Growth Conjecture}
Let $\mathbb{P}^n(\Q)$ denote the n-dimensional projective space over the rationals $\Q$. For a point $p=[\bfp]\in \mathbb{P}^n(\Q)$, with $\bfp\in \Z^{n+1}$ chosen so that $\gcd(p_0,p_1,\cdots,p_n)=1$, we define its height to be 
\begin{align}
    H(p):=\max_{0\leq i\leq n} |p_i| .
\end{align}
Consider an irreducible variety $V\subset \mathbb{P}^n(\Q)$ of degree $d\ge 1$, we define its counting function as 
\begin{align}
    N_V(B):=\# \{p\in V:H(p)\leq B \} .
\end{align}
In \cite{Ser16}, Serre proposed the following conjecture:
\begin{conjecture}[Serre's Dimension Growth Conjecture]
    Let $V\subset \mathbb{P}^n(\Q)$ be an irreducible projective variety of degree $d\ge 2$ defined over $\Q$. Then
    \begin{align}
        N_V(B)\lesim_{V} B^{\dim V}(\log B)^c
    \end{align}
    for some constant $c=c_V>0$.
\end{conjecture}
Our main result Theorem 1.2 has a direct corollary that can be regarded as the analogue of the dimension growth conjecture for smooth manifolds. Let $X$ be a compact submanifold of $\R^n$, and let $(\bfx,\mathbf{f}(\bfx))$ be a local chart of $X$, where $\mathbf{f}=(f_1,f_2,\cdots,f_k)$ and $\bfx=(x_1,x_2,\cdots,x_m)$. When $m\ge 2$, Srivastava and Technau \cite{ST23} proved Smooth variant of Serre’s Dimension Growth estimates for a family of flat manifold. In this paper, we focus on the case when $m=1$.

We say $X$ is of type $d\ge 2$ at $x_0$ if for every small neighbourhood $(x_0-r_0,x_0+r_0)$ with $r_0>0$ and any local chart $(x,\mathbf{f}(x))$ of $X$, where $x\in (x_0-r_0,x_0+r_0)$, there exist $\bfa \in \Z^{k} $ and $Q_0\in \N$ such that
\begin{align}
    x \to \mathbf{f}(x)\cdot \frac{\bfa}{Q_0}
\end{align}
is a curve of type $d\ge 2$ at $x_0$. We say $X$ is finite type if $X$ satisfies above condition everywhere. Then by letting $\delta=0$ in Theorem 1.2 and routine process(see for instance \cite{Hua20} Theorem 4 and \cite{ST23} Theorem 2.11), we have
\begin{theorem}
    Let $X$ be a compact submanifold of $\R^n$. Assume $X$ is finite type, then
    \begin{align}
        N_X(B) \lesim_{X,\epsilon} B^{\dim X+\epsilon}
    \end{align}
    for every $\epsilon>0$.
\end{theorem}

\section{Reduction to a model case}
By Taylor expansion, if $(t,f(t))$ is of finite type d at $a_0$ and $f\in C^{d+1}$, then there exists a $\epsilon_0>0$, such that
\begin{align}\label{2.1}
    f(t)=c_0+c_1(t-a_0)^d+R(t):=h(t)+R(t)
\end{align}
for all $t\in I_0:=(a_0-\epsilon_0,a_0+\epsilon_0)$, where $c_0, c_1$ are constants and $c_1\neq 0$. The remainder term is given by
\begin{align}
    R(t)=\frac{f^{(d+1)}(\xi)}{(d+1)!}(t-a_0)^{d+1},
\end{align}
where $\xi$ is between $a_0$ and $t$. By choosing $\epsilon_0$ sufficiently small, we can assume the following
\begin{enumerate}
    \item [(H1)]  
    \begin{align}\label{H1}       
    |f^{(d)}(t)|\gtrsim 1,\ \forall t\in I_0 .
    \end{align}
    \item [(H2)] For $j=1,\cdots,d$
    \begin{align}
        |f^{(j)}(t)|\simeq |h^{(j)}(t)|\simeq 1, \ \forall t\in (a_0-\epsilon_0,a_0-\frac{\epsilon_0}{2})\cup(a_0+\frac{\epsilon_0}{2},a_0+\epsilon_0) .
    \end{align}
    \item [(H3)]\label{H3} On each half interval $(a_0-\epsilon_0,a_0)$ and $(a_0,a_0+\epsilon_0)$, $f'$ is either strictly increasing or strictly decreasing.
\end{enumerate}
We refer to such $f:I_0 \to \R$ as the model case.

We will prove the following result for the model case
\begin{proposition}\label{Pro2.1}
    Let $f:I_0 \to \R$ defined in \eqref{2.1} satisfy (H1) to (H3), then for any $\delta \in (Q^{-1},\frac{1}{2})$ and integer $Q\ge 1$, we have 
    \begin{align}
          N_f(Q,\delta)=|I|\delta Q^2+O\Big{(}\delta^{\frac{1}{d}}Q^{2-\frac{1}{d}}+Q^{1+\epsilon }\Big{)}
    \end{align}
    for all $\epsilon>0$.
\end{proposition}

\section{Proof of the model case}
Now, we focus on the model case. Let $I_0=(a_0-\epsilon_0,a_0+\epsilon_0)$, and let $f:I_0\to \R$ be a $C^{d+1}$ function such that
\begin{align}
    f(t)=c_0+c_1(t-a_0)^d+\frac{f^{(d+1)}(\xi)}{(d+1)!}(t-a_0)^{d+1} ,
\end{align}
where $\xi$ is between $a_0$ and $t$, $c_1\neq 0$. Assuming $f$ is the model case, we prove the following result:
\begin{theorem}
   For any $\delta \in (Q^{-1},\frac{1}{2})$ and integer $Q\ge 1$, we have 
   \begin{align}     
    N_f(Q,\delta)=|I|\delta Q^2+O\Big{(}\delta^{\frac{1}{d}}Q^{2-\frac{1}{d}}+Q^{1+\epsilon }\Big{)}  .
   \end{align}
\end{theorem}

Let $\{u_n\}$ be a sequence of $N$ numbers with $1\leq n\leq N$. Let $Z(N,\alpha,\beta)$ count the number of $n$ for which $u_n \in (\alpha,\beta)$ (mod 1) with $\beta \in (\alpha,\alpha+1)$. Let 
\begin{align}
    D(N,\alpha,\beta):=Z(N,\alpha,\beta)-(\beta-\alpha)N
\end{align}
be the discrepancy function. We use the following lemma:
\begin{lemma}[Chapter 1 \cite{Mon94}]
    For any positive integer $K$, we have
    \begin{align}
        |D(N,\alpha,\beta)|\leq \frac{N}{K+1}+2\sum_{k=1}^K b_k \Big{|} \sum_{n=1}^N e(ku_n) \Big{|} ,
    \end{align}
    where 
\begin{align}
    b_k=\frac{1}{K+1}+\min \Big{(} \beta-\alpha,\frac{1}{\pi k} \Big{)}.
\end{align}
\end{lemma}

We apply Lemma 2.2 to the sequence $\{qf(\frac{a}{q})  \}$, where $\frac{a}{q}\in I_0$ and $q\leq Q$, with $\alpha=-\delta$ and $\beta=\delta$. The total number of such rational numbers is 
\begin{align}
    \sum_{q\leq Q}\sum_{\frac{a}{q}\in I_0} 1= \sum_{q\leq Q}(|I_0|q+O(1))= \frac{|I_0|Q^2}{2}+O(Q).
\end{align}
Therefore for any positive integer $K$, we have
\begin{align}
    N_f(Q,\delta)-|I_0|\delta Q^2 \lesim \frac{Q^2}{K}+\delta Q+\sum_{k=1}^K b_k \Big{|}\sum_{q\leq Q}\sum_{\frac{a}{q}\in I_0}e\Big(kqf(\frac{a}{q})\Big) \Big{|} ,
\end{align}
where 
\begin{align}
    b_k=\frac{1}{K+1}+\min \Big{(} 2\delta,\frac{1}{\pi k} \Big{)}.
\end{align}
We will choose
\begin{align}
    K=Q.
\end{align}
Since $f$ satisfies (H3), by dividing the interval $I_0$ into two half intervals $(a_0-\epsilon_0,a_0)$ and $(a_0,a_0+\epsilon_0)$,  we can assume $f'$ is either strictly increasing or decreasing. Without loss of generality, we only need to deal with case when $\frac{a}{q}\in (a_0,a_0+\epsilon_0)$. We still denote $I_0:=(a_0,a_0+\epsilon_0)$, we can further assume $f'>0$ and strictly increasing on $I$.

We use the following truncated version of Poisson summation formula.
\begin{lemma}[Chapter 3 \cite{Mon94}]
    Let $g$ be a real-valued function and suppose that $g'$ is continuous and increasing on $[c,d]$. Let $s=g'(c)$ and $t=g'(d)$. Then
    \begin{align}
        \sum_{c\leq n\leq d} e\Big(g(n)\Big)=\sum_{s-1\leq j\leq t+1}\int_{c}^d e\Big(g(x)-jx\Big)dx+O(\log (2+t-s)) .
    \end{align}
\end{lemma}
Denote $J:=[\inf_{t\in I_0} f'(t),\sup_{t\in I_0} f'(t)]$ and $J_k:=[k\inf_{t\in I_0} f'(t)-1,k\sup_{t\in I_0} f'(t)+1]$. Then by Lemma 2.3, we have
\begin{align}
    \sum_{a\in qI_0} e\Big(kqf(\frac{a}{q})\Big)=\sum_{j\in J_k} \int_{qI_0}e\Big(kqf(\frac{r}{q})-jr\Big)dr+O(\log(2+k|J|)).
\end{align}
Since $|J|\lesim |I_0|$, the contribution from the error term is
\begin{align}
    \lesim Q\sum_{k=1}^K\Big{(} \frac{1}{K+1}+\min \Big{(} 2\delta, \frac{1}{\pi k}\Big{)}  \Big{)}\log k \lesim Q(\log K)^2 .
\end{align}
By the change of variables $r=qx$, the inner integral of the main term becomes
\begin{align}
    q\int_{I_0} e\Big(q(kf(x)-jx)\Big)dx.
\end{align}

Let us deal with some endpoint cases. To avoid cases where there is a critical point very close to an endpoint of $I_0$, we shrink the interval $J_k$ slightly. Let $\Tilde{J}_k:=[k\inf_{t\in I_0} f'(t)+1,k\sup_{t\in I_0} f'(t)-1]$; if $k\sup_{t\in I_0} f'(t)-1<k\inf_{t\in I_0} f'(t)+1$, we take $\Tilde{J}_k$ to be empty. We need the following lemma for oscillatory integrals.
\begin{lemma}[Van der Corput's lemma]
   Let $k\ge 2$ be an integer. Suppose $\phi:(a,b)\to \R$ is $C^k$, and that
   \begin{align}
       |\phi^{(k)}(x)|\ge 1 
   \end{align}
   for all $x\in (a,b)$. Then 
   \begin{align}
       \Big| \int_{a}^b e^{i\lambda \phi(x)}dx \Big| \leq c_k \lambda^{-\frac{1}{k}}.
   \end{align}
   The constant $c_k$ is independent of $\phi$ and $\lambda$.
\end{lemma}
By van der Corput's lemma,
\begin{align}
    \int_{I_0} e\Big(q(kf(x)-jx)\Big)dx \lesim (qk)^{-\frac{1}{d}},
\end{align}
since 
\begin{align}
    f^{(d)}(x)\gtrsim 1, \forall x\in I_0.
\end{align}
Thus, we get
\begin{align}
    &\sum_{k\leq K} b_k \Big{|} \sum_{q\leq Q} \sum_{j\in J_k \backslash \Tilde{J}_k} q\int_{I_0} e\Big(q(kf(x)-jx)\Big)dx\Big{|}\\
    &\lesim \sum_{k\leq K}\Big{(}\frac{1}{K+1}+\min \Big{(} \delta,\frac{1}{k}\Big{)} \Big{)} k^{-\frac{1}{d}}\sum_{q\leq Q}q^{1-\frac{1}{d}}\\
    &\lesim K^{-\frac{1}{d}}Q^{2-\frac{1}{d}}+\delta^{\frac{1}{d}}Q^{2-\frac{1}{d}}\\
    &\lesim \delta^{\frac{1}{d}}Q^{2-\frac{1}{d}} ,
\end{align}
since $K\ge \frac{1}{\delta}$.

We now bound the following expression:
\begin{align}
    \sum_{k\leq K} b_k \Big{|} \sum_{q\leq Q}\sum_{j\in \Tilde{J}_k}q\int_{I_0} e\Big(qk(f(x)-\frac{j}{k}x)\Big)dx \Big{|} .
\end{align}
We perform a dyadic decomposition in $\frac{j}{k}$ and $x$. Let $\chi \in C_c^{\infty}(\R)$ be a bump function supported on $(-2\epsilon_0,2\epsilon_0)$ such that $\chi(t)=1$ for $t\in (-\epsilon_0,\epsilon_0)$. Define $\chi_k(t):=\chi(2^kt)-\chi(2^{k+1}t)$ for $k\ge 0$. Then the above summation is bounded by
\begin{align}
    \lesim I_1+I_2 ,
\end{align}
where 
\begin{align}
    &I_1:=\sum_{k\leq K} b_k \Big{|} \sum_{m=0}^M \sum_{n=0}^{N} \sum_{q\leq Q}\sum_{j\in \Tilde{J}_{k,m}}q\int_{I_0} e\Big(qk(f(x)-\frac{j}{k}x)\Big)\chi_n(x-a_0) dx \Big{|} , \\
    &I_2:=\sum_{k\leq K} b_k \Big{|} \sum_{m=0}^M \sum_{q\leq Q}\sum_{j\in \Tilde{J}_{k,m}}q\int_{I_0} e\Big(qk(f(x)-\frac{j}{k}x)\Big)\chi(2^{N+1}(x-a_0)) dx \Big{|},
\end{align}
where $M \simeq \log_2 k \leq \log_2 K $, $\Tilde{J}_{k,m}:=\{j\in \Tilde{J}_k: \frac{j}{k}\simeq 2^{-m} \}$ and $N:=\frac{m}{d-1}$. 

We do a change of variable $x\to x+a_0$, denote $I:=(0,\epsilon_0)$, then 
\begin{align}
    &I_1:=\sum_{k\leq K} b_k \Big{|} \sum_{m=0}^M \sum_{n=0}^{N} \sum_{q\leq Q}\sum_{j\in \Tilde{J}_{k,m}}q\int_{I} e\Big(qk(f(x+a_0)-\frac{j}{k}(x+a_0))\Big)\chi_n(x) dx \Big{|} , \\
    &I_2:=\sum_{k\leq K} b_k \Big{|} \sum_{m=0}^M \sum_{q\leq Q}\sum_{j\in \Tilde{J}_{k,m}}q\int_{I} e\Big(qk(f(x+a_0)-\frac{j}{k}(x+a_0))\Big)\chi(2^{N+1}x) dx \Big{|} .
\end{align}

We recall two basic lemmas for oscillatory integrals, which can be found in \cite{Ste93} or \cite{Hua20}.
\begin{lemma}[Non-stationary phase]
    Let $u\in C_c^{\infty}(\R^d)$ and $\phi \in C^l(\R^d)$ with $\nabla \phi(\bfx) \neq 0$ for $\bfx \in \supp(u)$. Then for any $\lambda>0$, we have
    \begin{align}
        \Big| \int_{\R^d} u(\bfx) e(\lambda \phi(\bfx))d\bfx \Big| \leq C_l \lambda^{-l+1},
    \end{align}
    and $C_l$ depends only on upper bounds for finitely many derivatives of $u$ and $\phi$ and a lower bound for $\nabla \phi$ in the support of $u$.
\end{lemma}

\begin{lemma}[Stationary phase]
    Let $u\in C_c^{\infty}(\R^d)$ and $\phi \in C^l(\R^d)$ for some integer $l>\frac{d}{2}+4$. Suppose that $\nabla \phi(\bfx_0)=\bf0$ and $\det \nabla^2 \phi(\bfx_0)\neq 0$. Let $\sigma$ be the signature of $\nabla^2 \phi(\bfx_0)$ and $\Delta=|\det \nabla^2 \phi(\bfx_0)|$. Suppose further $\nabla \phi(\bfx)\neq \bf0$ for $\bfx \in \supp (u)\backslash \{\bfx_0\} $. Then for $\lambda>0$, we have
    \begin{align}
        \int_{\R^d} u(\bfx) e(\lambda \phi(\bfx))d\bfx= e\Big(\lambda \phi(\bfx_0)+\frac{\sigma}{8}\Big) \Delta^{-\frac{1}{2}} \lambda^{-\frac{d}{2}} \Big(u(\bfx_0)+O(\lambda^{-1})\Big) ,
    \end{align}
    where the implicit constant depends only on upper bounds for finitely many derivatives of $u$ and $\phi$ and a lower bound for $\Delta$.
\end{lemma}

For $I_2$, since $f$ is a model case, the first order derivative of the phase function never vanishes,
\begin{align}
    |qk(f(x+a_0)-\frac{j}{k}x)'| \gtrsim qj .
\end{align}
Thus we can use integration by parts to finish the proof. We now focus on $I_1$. Let 
\begin{align}
    f_{n,jk}(x):= 2^{nd}f(2^{-n}x+a_0)-\frac{j}{k}2^{nd}a_0:=f_n(x)-\frac{j}{k}2^{nd}a_0.
\end{align}
Since $f_{n}$ is a model case, we have
\begin{align}
    I_1\lesim \sum_{k\leq K} b_k \Big{|} \sum_{m=0}^M \sum_{n=0}^N \sum_{q\leq Q}\sum_{j\in \Tilde{J}_{k,m}}q2^{-n}\int_I e\Big(qk2^{-nd}(f_{n,jk}(x)-\frac{j}{k}2^{n(d-1)}x)\Big)\chi_0(x) dx \Big{|} .
\end{align}
By the uncertainty principle, it's natural to separate the summation in $q$ into 2 parts. Denote
\begin{align}
    I_{1,1}:= \sum_{k\leq K} b_k \Big{|} \sum_{m=0}^M \sum_{n=0}^N \sum_{2^{nd}k^{-1} \leq q\leq Q}\sum_{j\in \Tilde{J}_{k,m}}q2^{-n}\int_I e\Big(qk2^{-nd}(f_{n,jk}(x)-\frac{j}{k}2^{n(d-1)}x)\Big)\chi_0(x) dx \Big{|}
\end{align}
and
\begin{align}
    I_{1,2} := \sum_{k\leq K} b_k \Big{|} \sum_{m=0}^M \sum_{n=0}^N \sum_{1 
\leq q\leq 2^{nd}k^{-1}}\sum_{j\in \Tilde{J}_{k,m}}q2^{-n}\int_I e\Big(qk2^{-nd}(f_{n,jk}(x)-\frac{j}{k}2^{n(d-1)}x)\Big)\chi_0(x) dx \Big{|} .
\end{align}
For $I_{1,2}$, since there is no oscillation in the phase function, we apply the trivial triangle inequality
\begin{align}
    I_{1,2} &\lesim \sum_{k\leq K} b_k \Big{|}   
  \sum_{m=0}^M \sum_{j\in \Tilde{J}_{k,m}} \sum_{n=0}^N \sum_{1\leq q \leq 2^{nd}k^{-1}} q2^{-n}  \Big{|} \\
    & \lesim \sum_{k\leq K} b_k \Big{|}   
  \sum_{m=0}^M \sum_{j\in \Tilde{J}_{k,m}} (k^{-1})^2 2^{N(2d-1)} \Big{|}\\
  & \lesim \sum_{k\leq K} b_k |k^{-1}\sum_{m=0}^M 2^{m\frac{d}{d-1}}|\\
  & \lesim \sum_{k\leq K} b_k k^{\frac{1}{d-1}} \lesim K^{\frac{1}{d-1}} ,
\end{align}
which is better than what we need. Therefore we only need to focus on $I_{1,1}$. Notice that
\begin{align}
    \frac{j}{k}2^{n(d-1)}\leq 1 .
\end{align}
Since $f_{n}'$ is increasing, we denote the unique solution of 
\begin{align}
    f_{n}'(x)=\frac{j}{k}2^{n(d-1)}
\end{align}
by $\theta_1$. Notice that we only need to consider the case $n=N$, as the cut-off function is supported on $x\simeq 1$. Since $f_{N}''\gtrsim 1$ on the support of $\chi_0$, by the stationary phase method we have
\begin{align}
    &I_{1,1} \lesim\\
    &\sum_{k\leq K} b_k \Big{|} \sum_{m=0}^M \sum_{2^{Nd}k^{-1} \leq q\leq Q}\sum_{j\in \Tilde{J}_{k,m}}q^{\frac{1}{2}}k^{-\frac{1}{2}}2^{\frac{Nd}{2}-N}   e\Big(qk2^{-Nd}(f_{N,jk}(\theta_1)-\frac{j}{k}2^{N(d-1)}\theta_1)\Big)a_{N}(\theta_1) \Big{|}\\
    &+ error ,
\end{align}
where $a_{N}(t):=|f_{N}''(t)|^{-\frac{1}{2}}\chi_0(t)$. We denote the first term by $I_{main}$. Let's deal with the error term:
\begin{align}
    error &\lesim \sum_{k\leq K} b_k \sum_{m=0}^M \sum_{2^{Nd}k^{-1} \leq q\leq Q}\sum_{j\in \Tilde{J}_{k,m}}q2^{-N}(qk2^{-Nd})^{-1}\\
    &\lesim Q\log K\sum_{k\leq K} b_k .
\end{align}
Recall that
\begin{align}
    b_k=\frac{1}{K+1}+\min \Big{(} 2\delta, \frac{1}{\pi k} \Big{)} \lesim \frac{1}{K+1}+\frac{1}{k}.
\end{align}
As a consequence,
\begin{align}
    error \lesim Q(\log Q)^2 .
\end{align}

Define the dual function $f^*$ of $f$ to be
\begin{align}
    f^*(y)=y(f')^{-1}(y)-f((f')^{-1}(y)) .
\end{align}
Notice that 
\begin{align}
    f_{N}(\theta_1)-\frac{j}{k}2^{N(d-1)}\theta_1=-f_{N}^*(\frac{j}{k}2^{N(d-1)}) .
\end{align}
We have
\begin{align}
    &I_{main} \lesim \\ 
    &\sum_{k\leq K} b_k \cdot k^{-\frac{1}{2}} \sum_{m=0}^M \Big{|} \sum_{2^{Nd}k^{-1} \leq  q\leq Q}\sum_{j\in \Tilde{J}_{k,m}} q^{\frac{1}{2}} 2^{-N+\frac{Nd}{2}} e\Big(-qk2^{-Nd}(f_{N}^*(\frac{j}{k}2^{N(d-1)})+\frac{j}{k}2^{Nd}a_0 )\Big)a_{N}(\theta_1) \Big{|}.
\end{align}
Denote 
\begin{align}
   & F(t):=f_{N}^*(t)+2^{N}a_0 t,\\ 
   & \delta^*:= \frac{2^N}{Q} \leq 1.
\end{align}
We first estimate the summation in $q$, we have
\begin{align}
   & \sum_{2^{Nd}k^{-1} \leq  q\leq Q} q^{\frac{1}{2}} e\Big(-qk2^{-m-N}F(\frac{j}{k}2^m)\Big) \\
    &\lesim
        \left\{\begin{array}{lr}
        Q^{\frac{3}{2}}  & \text{if}\ ||2^{-m}k F(\frac{j}{k}2^m) || \leq \delta^*, \\
        Q^{\frac{1}{2}}||2^{-m}kF(\frac{j}{k}2^m)||^{-1} & \text{if}\ ||2^{-m}kF(\frac{j}{k}2^m)|| >\delta^*.    \end{array}
    \right.        
\end{align}

This problem reduces to a rescaled version of counting rational points near curves with $F'' \simeq 1$, and thus we can apply the results of \cite{Hua15} or \cite{Hux94}. By a dyadic decomposition in $k$, we have 
\begin{align}
    I_{main}\lesim  \sum_{K_0} \Big( \Sigma_{1,K_0}+\Sigma_{2,K_0} \Big) ,
\end{align}
where 
\begin{align}
   & \Sigma_{1,K_0}:= Q^{\frac{3}{2}} \sum_{m=0}^M 2^{-N+\frac{Nd}{2}} \sum_{k\simeq K_0} b_k \cdot k^{-\frac{1}{2}}\sum_{\underset{||2^{-m}k F(\frac{j}{k}2^m) || \leq \delta^*}{j\in \Tilde{J}_{k,m}}} 1,  \\
   &\Sigma_{2,K_0}:=Q^{\frac{1}{2}} \sum_{m=0}^M 2^{-N+\frac{Nd}{2}} \sum_{k\simeq K_0} b_k \cdot k^{-\frac{1}{2}}\sum_{\underset{||2^{-m}k F(\frac{j}{k}2^m) || > \delta^* }{j\in \Tilde{J}_{k,m}}}||2^{-m}kF(\frac{j}{k}2^m)||^{-1} ,
\end{align}
and the summation is over all dyadic $K_0$ with $1\leq K_0\leq K$. We bound the first term $\Sigma_{1,K_0}$.

\begin{claim}\label{Claim4.7}
    We have 
    \begin{align}\label{M1}
        \sum_{k\simeq K_0} \sum_{\underset{||2^{-m}k F(\frac{j}{k}2^m) || \leq \delta^*}{j\in \Tilde{J}_{k,m}}} 1 \lesim_{\epsilon} 2^{-m}\delta^* K_0^2+K_0^{1+\epsilon}.
    \end{align}
\end{claim}

\begin{proof}
When $2^{-m} K_0\leq 1 $, the claim is trivial, since there is essentially one term in the summation over $j$. So without loss of generality, we assume $2^{-m} K_0>1 $. Let $\chi_\delta(\theta)$ be the characteristic function of the set 
\begin{align}
    \{ \theta\in \R: ||\theta||\leq \delta  \}
\end{align}
and 
\begin{align}
    J=\Big \lfloor \frac{1}{2\delta} \Big\rfloor.
\end{align}
Then consider the Fejér kernel
\begin{align}
    \mathcal{F}_J (\theta):= \frac{1}{J^2} \Big| \sum_{j=1}^J e(j\theta)   \Big|^2=\Big( \frac{\sin (\pi J\theta)}{J \sin(\pi \theta)} \Big)^2=\sum_{j=-J}^J \frac{J-|j|}{J^2} e(j\theta).
\end{align}
An easy computation shows that when $||\theta||\leq \delta$, we have
\begin{align}
    \Big( \frac{\sin (\pi J\theta)}{J \sin(\pi \theta)} \Big)^2 \ge \Big( \frac{2\pi^{-1}\pi J||\theta||}{J \pi ||\theta||} \Big)^2 =\frac{4}{\pi^2},
\end{align}
so 
\begin{align}
    \chi_\delta (\theta)\leq \frac{\pi^2}{4} \mathcal{F}_J(\theta).
\end{align}
Let $\omega$ be a smooth bump function adapted to $x\simeq 1$, then left hand side of $\eqref{M1}$ is 
\begin{align}
   & \lesim \sum_{k\simeq K_0} \sum_{\underset{||2^{-m}k F(\frac{j}{k}2^m) || \leq \delta^*}{j \in \Z }} w(\frac{j}{k}2^m) \\
   & \lesim \sum_{k\simeq K_0} \sum_{j \in \Z } w(\frac{j}{k}2^m) \sum_{j_0=-J_0}^{J_0} \frac{J_0-|j_0|}{J_0^2} e\Big( j_0k 2^{-m} F(\frac{j}{k}2^m)  \Big),
\end{align}
where $J_0:= \lfloor \frac{1}{2\delta^*} \rfloor $. When $j_0=0$, right hand side contributes
\begin{align}
    \simeq   2^{-m}\delta^* K_0^2.
\end{align}
We only need to bound 
\begin{align}
   \sum_{j_0\simeq J_1} \frac{J_0-j_0}{J_0^2} \sum_{k\simeq K_0} \sum_{j \in \Z } w(\frac{j}{k}2^m)  e\Big( j_0k 2^{-m} F(\frac{j}{k}2^m)  \Big)
\end{align} 
for dyadic $J_1$ with $1\leq J_1\leq J_0$. By Poisson summation, we have
\begin{align}
     \sum_{j \in \Z } w(\frac{j}{k}2^m)  e\Big( j_0k 2^{-m} F(\frac{j}{k}2^m)  \Big) &=\sum_{a\in \Z} \int_{\R} w(\frac{t}{k}2^m) e\Big( j_0k 2^{-m} F(\frac{t}{k}2^m)-at \Big)  dt \\
     &=k 2^{-m} \sum_{a\in \Z} \int_{\R} w(x) e\Big( j_0k 2^{-m} (F(x)-\frac{a}{j_0}x) \Big)  dx.
\end{align}
We take a derivative of the phase function, we have
\begin{align}
    \Big(F(x)-\frac{a}{j_0} x  \Big)'=(f_N')^{-1}(x)+2^N a_0-\frac{a}{j_0},
\end{align}
where we use
\begin{align}
    (f_N^*)'(x)=(f_N')^{-1}(x).
\end{align}
By non-stationary phase principle, we only need to focus on those $a$ with 
\begin{align}
   \Big| \frac{a}{j_0}-2^N a_0 \Big| \simeq 1.
\end{align}
For these $a$, by stationary phase, we have
\begin{align}
    \int_{\R} w(x) e\Big( j_0k 2^{-m} (F(x)-\frac{a}{j_0}x) \Big)\simeq (j_0k 2^{-m})^{-\frac{1}{2}} e\Big(-j_0k 2^{-m} F^*(\frac{a}{j_0}) \Big)+\text{error}.
\end{align}
The error terms can be bounded easily. The main term is 
\begin{align}\label{Main_term}
   2^{-\frac{m}{2}}  \sum_{j_0\simeq J_1} \frac{J_0-j_0}{J_0^2} j_0^{-\frac{1}{2}}\sum_{|\frac{a} {j_0}-2^N a_0|\simeq 1}  \sum_{k\simeq K_0} k^{\frac{1}{2}} e\Big(-j_0k 2^{-m} F^*(\frac{a}{j_0}) \Big).
\end{align}
We use 
\begin{align}
    \sum_{k\simeq K_0} k^{\frac{1}{2}} e\Big(-j_0k 2^{-m} F^*(\frac{a}{j_0}) \Big) \lesim
        \left\{\begin{array}{lr}
        K_0^{\frac{3}{2}}  & \text{if}\ ||j_0 F^*(\frac{a}{j_0}) || \leq (2^{-m}K_0)^{-1}, \\
        K_0^{\frac{1}{2}}||j_0 F^*(\frac{a}{j_0})||^{-1} & \text{if}\ ||j_0 F^*(\frac{a}{j_0}))|| >(2^{-m}K_0)^{-1} .    \end{array}
    \right. 
\end{align}
Then one would conclude that
\begin{align}
    \eqref{Main_term} \lesim 2^{-\frac{m}{2}}\delta^*J_1^{-\frac{1}{2}}\Big(K_0^{\frac{3}{2}} \sum_{j_0\simeq J_1} \sum_{\underset{||j_0 F^*(\frac{a}{j_0}) || \leq (2^{-m}K_0)^{-1} }{|\frac{a} {j_0}-2^N a_0|\simeq 1}} 1  +  K_0^{\frac{1}2{}} \sum_{j_0\simeq J_1} \sum_{\underset{||j_0 F^*(\frac{a}{j_0}) || > (2^{-m}K_0)^{-1} }{|\frac{a} {j_0}-2^N a_0|\simeq 1}}  ||j_0 F^*(\frac{a}{j_0}) ||^{-1}  \Big).
\end{align}
We bound the first term. Notice that 
\begin{align}\label{parabola}
    \sum_{j_0\simeq J_1} \sum_{\underset{||j_0 F^*(\frac{a}{j_0}) || \leq (2^{-m}K_0)^{-1} }{|\frac{a} {j_0}-2^N a_0|\simeq 1}} 1
\end{align}
is counting rational point near a convex curve, which can be evaluated by the main result of \cite{VV06}.
\begin{theorem}[\cite{VV06}, \cite{Hua15}]
    Suppose that $f\in C^3(I)$ and has a second
derivative bounded away from 0 on a bounded interval $I$. Let $\delta \in (0,\frac{1}{2})$, then for any $\epsilon>0$ and $Q\ge 1$,
\begin{align}
    \sum_{q\simeq Q} \sum_{\underset{||qf(\frac{a}{q})||\leq \delta}{\frac{a}{q}\in I}} 1 \lesim_{\epsilon} \delta Q^2+ Q^{1+\epsilon}.
\end{align}
\end{theorem}
We need to check derivative conditions. Since 
\begin{align}
   & (F^*)'(x)= (F')^{-1}(x)=f_N'(x-2^N a_0) , \\
   & (F^*)''(x)= f_N''(x-2^Na_0),
\end{align}
we get
\begin{align}
    (F^*)'(\frac{a}{j_0}),\ (F^*)''(\frac{a}{j_0}) \simeq 1
\end{align}
whenever $|\frac{a} {j_0}-2^N a_0|\simeq 1$. So we have
\begin{align}
    \eqref{parabola} \lesim_{\epsilon} (2^{-m}K_0)^{-1}J_1^2  +J_1^{1+\epsilon}.
\end{align}
And thus 
\begin{align}\label{M}
    2^{-\frac{m}{2}}\delta^*J_1^{-\frac{1}{2}} K_0^{\frac{3}{2}} \sum_{j_0\simeq J_1} \sum_{\underset{||j_0 F^*(\frac{a}{j_0}) || \leq (2^{-m}K_0)^{-1} }{|\frac{a} {j_0}-2^N a_0|\simeq 1}} 1  \lesim_{\epsilon}  2^{-\frac{m}{2}}\delta^*J_1^{-\frac{1}{2}} K_0^{\frac{3}{2}} \Big( (2^{-m}K_0)^{-1}J_1^2  +J_1^{1+\epsilon} \Big) .
\end{align}

The second term can be dealt with the same computation and a dyadic decomposition in $||j_0 F^*(\frac{a}{j_0}) ||$. We would get 
\begin{align}
    &2^{-\frac{m}{2}}\delta^*J_1^{-\frac{1}{2}} K_0^{\frac{1}{2}} \sum_{j_0\simeq J_1} \sum_{\underset{||j_0 F^*(\frac{a}{j_0}) || > (2^{-m}K_0)^{-1} }{|\frac{a} {j_0}-2^N a_0|\simeq 1}}  ||j_0 F^*(\frac{a}{j_0}) ||^{-1} \\
    &\lesim_\epsilon 2^{-\frac{m}{2}} \delta^*J_1^{-\frac{1}{2}} K_0^{\frac{1}{2}} (2^{-m}K_0) \Big( (2^{-m}K_0)^{-1}J_1^2  +J_1^{1+\epsilon}  \Big) \lesim \ \text{right hand side of\ } \eqref{M} .
\end{align}
We put this back to \eqref{Main_term}, then sum in dyadic $J_1\in [1,J_0]$. We have
\begin{align}
   \lesim_{\epsilon} 2^{\frac{m}{2}}K_0^{\frac{1}{2}} (\delta^*)^{-\frac{1}{2}}+ 2^{-\frac{m}{2}} (\delta^*)^{\frac{1}{2}-\epsilon}K_0^{\frac{3}{2}}.
\end{align}
When $\delta^* \ge 2^m K_0^{-1}$, the above is bounded by
\begin{align}
    \lesim_{\epsilon} 2^{-\frac{m}{2}} (\delta^*)^{\frac{1}{2}-\epsilon}K_0^{\frac{3}{2}}+K_0,
\end{align}
so the whole counting is
\begin{align}
     \sum_{k\simeq K_0} \sum_{\underset{||2^{-m}k F(\frac{j}{k}2^m) || \leq \delta^*}{j\in \Tilde{J}_{k,m}}} 1 \lesim_{\epsilon} 2^{-m}\delta^* K_0^2+K_0^{1+\epsilon}. \
\end{align}
When $\delta^*< 2^mK_0^{-1}$, we use the monotonicity property of the counting problem in $\delta^*$ to finish the proof.
    
\end{proof}

Now we return to estimate $\Sigma_{1,K_0}$. Recall that 
\begin{align}
    & b_k=\frac{1}{K+1}+\min \Big{(} 2\delta,\frac{1}{\pi k} \Big{)} \lesim \left\{\begin{array}{lr}
       \delta  & \text{if}\ k \leq \delta^{-1} , \\
        \frac{1}{k} & \text{if}\ k \ge \delta^{-1} .    \end{array}
    \right. \\
    & 2^M \simeq K_0.
\end{align}
By applying Claim \ref{Claim4.7}, we get
\begin{align}
    \Sigma_{1,K_0} & \lesim K_0^{-\frac{3}{2}}Q^{\frac{3}{2}} \sum_{m=0}^M 2^{-N+\frac{Nd}{2}} \sum_{k\simeq K_0}  \sum_{\underset{||2^{-m}k F(\frac{j}{k}2^m) || \leq \delta^*}{j\in \Tilde{J}_{k,m}}} 1    \\
   & \lesim_{\epsilon} K_0^{-\frac{3}{2}}Q^{\frac{3}{2}} \sum_{m=0}^M 2^{-N+\frac{Nd}{2}}\Big(2^{-m}\delta^* K_0^2+K_0^{1+\epsilon} \Big)     \\
   & \lesim_\epsilon  K_0^{\frac{1}{2}}Q^{\frac{1}{2}}+ K_0^{-\frac{1}{2(d-1)}+\epsilon}Q^{\frac{3}{2}} 
\end{align}
whenever $K_0 \ge \delta^{-1}$ and 
\begin{align}
    \Sigma_{1,K_0} & \lesim \delta K_0^{-\frac{1}{2}}Q^{\frac{3}{2}} \sum_{m=0}^M 2^{-N+\frac{Nd}{2}} \sum_{k\simeq K_0}  \sum_{\underset{||2^{-m}k F(\frac{j}{k}2^m) || \leq \delta^*}{j\in \Tilde{J}_{k,m}}} 1  \\
    & \lesim_{\epsilon} \delta K_0^{\frac{3}{2}}Q^{\frac{1}{2}} + \delta K_0^{1-\frac{1}{2(d-1)}+\epsilon} Q^{\frac{3}{2}}
\end{align}
whenever $K_0 \leq \delta^{-1}$. As a result,
\begin{align}
    \sum_{K_0} \Sigma_{1,K_0}= \sum_{\delta^{-1}\leq K_0\leq K}\Sigma_{1,K_0} +\sum_{K_0 \leq \delta^{-1}} \Sigma_{1,K_0} \lesim_{\epsilon} Q^{1+\epsilon}+\delta^{\frac{1}{2(d-1)}-\epsilon}Q^{\frac{3}{2}} .
\end{align}
where the summation of $K_0$ is over all dyadic numbers between $1$ and $K$. We can bound the summation of $\Sigma_{2,K_0}$ in the same way
\begin{align}
    \sum_{K_0}\Sigma_{2,K_0} \lesim_{\epsilon} Q^{1+\epsilon}+\delta^{\frac{1}{2(d-1)}-\epsilon}Q^{\frac{3}{2}}.
\end{align}
Notice that 
\begin{align}
    \delta^{\frac{1}{2(d-1)}} Q^{\frac{3}{2}}
\end{align}
is a weighted geometric mean of $\delta^{\frac{1}{d}}Q^{2-\frac{1}{d}}$ and $Q$. We conclude that
\begin{align}
     N_f(Q,\delta)-|I_0|\delta Q^2 \lesim_{\epsilon} \delta^{\frac{1}{d}}Q^{2-\frac{1}{d}}+Q^{1+\epsilon}
\end{align}
for $\delta \in (Q^{-1},\frac{1}{2})$.

\normalem

\noindent Department of Mathematics, University of Wisconsin-Madison, Madison, WI-53706, USA \\
Email addresses:\\
mchen454@math.wisc.edu

\end{document}